\documentclass[11pt,fleqn,letter]{article}
\usepackage{amsfonts, epsfig}
\usepackage{color}
\newcommand{\ol}{\setlength{\itemsep}{0pt.}\begin{enumerate}}
\newcommand{\eol}{\end{enumerate}\setlength{\itemsep}{-\parsep}}
\newcommand{\ignore}[1]{}
\title{Bounds on the permanent and some applications}
\author{Leonid Gurvits\thanks{Department of Computer Science, Grove School of Engineering,
The City College of New York,
New York, United States of America. Research supported by NSF grant 116143.} ~~and~ Alex Samorodnitsky\thanks{
School of Engineering and Computer Science,
The Hebrew University of Jerusalem,
Jerusalem, Israel. Research supported by BSF and ISF grants.}
}

\begin{document}
\date{}
\maketitle


\newtheorem{THEOREM}{Theorem}[section]
\newenvironment{theorem}{\begin{THEOREM} \hspace{-.85em} {\bf :}
}%
                        {\end{THEOREM}}
\newtheorem{LEMMA}[THEOREM]{Lemma}
\newenvironment{lemma}{\begin{LEMMA} \hspace{-.85em} {\bf :} }%
                      {\end{LEMMA}}
\newtheorem{COROLLARY}[THEOREM]{Corollary}
\newenvironment{corollary}{\begin{COROLLARY} \hspace{-.85em} {\bf
:} }%
                          {\end{COROLLARY}}
\newtheorem{PROPOSITION}[THEOREM]{Proposition}
\newenvironment{proposition}{\begin{PROPOSITION} \hspace{-.85em}
{\bf :} }%
                            {\end{PROPOSITION}}
\newtheorem{DEFINITION}[THEOREM]{Definition}
\newenvironment{definition}{\begin{DEFINITION} \hspace{-.85em} {\bf
:} \rm}%
                            {\end{DEFINITION}}
\newtheorem{EXAMPLE}[THEOREM]{Example}
\newenvironment{example}{\begin{EXAMPLE} \hspace{-.85em} {\bf :}
\rm}%
                            {\end{EXAMPLE}}
\newtheorem{CONJECTURE}[THEOREM]{Conjecture}
\newenvironment{conjecture}{\begin{CONJECTURE} \hspace{-.85em}
{\bf :} \rm}%
                            {\end{CONJECTURE}}
\newtheorem{MAINCONJECTURE}[THEOREM]{Main Conjecture}
\newenvironment{mainconjecture}{\begin{MAINCONJECTURE} \hspace{-.85em}
{\bf :} \rm}%
                            {\end{MAINCONJECTURE}}
\newtheorem{PROBLEM}[THEOREM]{Problem}
\newenvironment{problem}{\begin{PROBLEM} \hspace{-.85em} {\bf :}
\rm}%
                            {\end{PROBLEM}}
\newtheorem{QUESTION}[THEOREM]{Question}
\newenvironment{question}{\begin{QUESTION} \hspace{-.85em} {\bf :}
\rm}%
                            {\end{QUESTION}}
\newtheorem{REMARK}[THEOREM]{Remark}
\newenvironment{remark}{\begin{REMARK} \hspace{-.85em} {\bf :}
\rm}%
                            {\end{REMARK}}

\newcommand{\thm}{\begin{theorem}}
\newcommand{\lem}{\begin{lemma}}
\newcommand{\pro}{\begin{proposition}}
\newcommand{\dfn}{\begin{definition}}
\newcommand{\rem}{\begin{remark}}
\newcommand{\xam}{\begin{example}}
\newcommand{\cnj}{\begin{conjecture}}
\newcommand{\mcnj}{\begin{mainconjecture}}
\newcommand{\prb}{\begin{problem}}
\newcommand{\fac}{\begin{fact}}

\newcommand{\que}{\begin{question}}
\newcommand{\cor}{\begin{corollary}}
\newcommand{\prf}{\noindent{\bf Proof:} }
\newcommand{\ethm}{\end{theorem}}
\newcommand{\elem}{\end{lemma}}
\newcommand{\epro}{\end{proposition}}
\newcommand{\edfn}{\bbox\end{definition}}
\newcommand{\erem}{\bbox\end{remark}}
\newcommand{\exam}{\bbox\end{example}}
\newcommand{\ecnj}{\bbox\end{conjecture}}
\newcommand{\emcnj}{\bbox\end{mainconjecture}}
\newcommand{\eprb}{\bbox\end{problem}}
\newcommand{\eque}{\bbox\end{question}}
\newcommand{\ecor}{\end{corollary}}
\newcommand{\eprf}{\bbox}
\newcommand{\efac}{\end{fact}}
\newcommand{\beqn}{\begin{equation}}
\newcommand{\eeqn}{\end{equation}}
\newcommand{\wbox}{\mbox{$\sqcap$\llap{$\sqcup$}}}
\newcommand{\bbox}{\vrule height7pt width4pt depth1pt}
\newcommand{\qed}{\bbox}
\def\sup{^}

\def\H{\{0,1\}^n}

\def\S{S(n,w)}

\def\g{\phi_{\ast}}
\def\y{y_{\ast}}
\def\z{z_{\ast}}
\def\x{x_{\ast}}

\def\f{\tilde f}

\def\n{\lfloor \frac n2 \rfloor}

\def \E{\mathbb E}
\def \R{\mathbb R}
\def \Z{\mathbb Z}
\def \N{\mathbb N}
\def \F{\mathbb F}
\def \T{\mathbb T}

\def \r{\textcolor{red}{r}}
\def \Rc{\textcolor{red}{R}}

\def \noi{{\noindent}}

\def \iff{~~~~\Leftrightarrow~~~~}

\def\<{\left<}
\def\>{\right>}
\def \({\left(}
\def \){\right)}

\def \e{\epsilon}
\def \l{\lambda}


\def\Tp{Tchebyshef polynomial}
\def\Tps{TchebysDeto be the maximafine $A(n,d)$ l size of a code with distance $d$hef polynomials}
\newcommand{\rarrow}{\rightarrow}

\newcommand{\larrow}{\leftarrow}

\overfullrule=0pt
\def\setof#1{\lbrace #1 \rbrace}

\begin{abstract}
We show that the permanent of a doubly stochastic $n \times n$ matrix
$A = \(a_{ij}\)$ is at least as large as $\prod_{i,j} \(1-a_{ij}\)^{1-a_{ij}}$
and at most as large as $2^n$ times this number. Combined with previous
work, this improves on the deterministic approximation factor for the
permanent, giving $2^n$ instead of $e^n$-approximation.

\noi We also give a combinatorial application of the lower bound, proving S. Friedland's "Asymptotic Lower
Matching Conjecture" for the monomer-dimer problem.
\end{abstract}

\section{Introduction}
The permanent of an $n \times n$ matrix $A = \(a_{ij}\)$ is given by
\vspace*{-2ex}
\[
Per(A) = \sum_{\sigma \in S_n} \prod_{i=1}^n a_{i\sigma(i)}
\]
\vspace*{-2ex} Here $S_n$ is the symmetric group on $n$ elements.

\noi The permanent is a classical mathematical notion, going back to Binet and Cauchy \cite{Minc}. One part of its appeal is its strong, though seemingly spurious, similarity to the determinant. Another part is in its ability to count things. The permanent of a $0,1$ matrix $A$ equals the number of perfect matchings in the bipartite graph it represents. The permanents are also useful in counting more complex subgraphs, such as Hamiltonian cycles (\cite{FKS} and the references therein).

\noi In fact, the permanent counts things in a very strong sense, since it is $\#P$ to compute \cite{Valiant}, even for $0,1$ matrices. Hence, from the complexity point of view, the permanent is very different from the determinant. While the latter is efficiently computable, the permanent of nonnegative matrices is (probably) not. The natural question is, therefore, to try and approximate the permanent as efficiently as possible, and as well as possible.

\noi We briefly discuss three different approaches to achieve this goal.

\noi {\bf The Monte Carlo Markov Chain approach}: As observed by Jerrum et al \cite{JSV} an efficient procedure to sample uniformly from the set of all perfect matchings in a bipartite graph is computationally equivalent to approximately counting the matchings. Broder \cite{Broder} proposed to construct such a procedure by devising a random walk on an appropriate space, rapidly converging to its stationary distribution, which would be uniform on the set of perfect matchings (and assign a substantial weight to it). This was accomplished (and extended) in \cite{JSV}, giving an efficient randomized approximation algorithm for the permanent of a nonnegative matrix, up to any degree of precision, and providing a complete solution to the problem.

\noi {\bf Exploiting the similarity to determinant}: This is based on an observation of Godsil and Gutman \cite{Lovasz-Plummer}, that, for a matrix $A = \(a_{ij}\)$ with nonnegative entries, the random matrix $B = \(\epsilon_{ij} \cdot \sqrt{a_{ij}}\)$ where $\epsilon_{ij}$ are independent random variables with expectation $0$ and variance $1$, satisfies $Per(A) = \E~Det^2(B)$.
Hence, for an efficient randomized permanent approximation, it would suffice to show the random variable $Det^2(B)$ to be concentrated around its expectation. In \cite{Barvinok} the random variables $\epsilon_{ij}$ were taken to be quaternionic Gaussians, leading to an efficient randomized approximation algorithm for the permanent, which achieves an approximation factor of about $1.3^n$.

\noi {\bf Using combinatorial bounds on the permanent}: The permanent of a doubly stochastic matrix was shown to be at least $\frac{n!}{n^n} \approx e^{-n}$ in \cite{Eg,Fal}, answering a question of van der Waerden. On the other hand, this permanent is (clearly) at most $1$. Hence, we already know the permanent of a doubly stochastic matrix up to a factor of $e^n$. In \cite{LSW} this fortuitous fact was exploited by showing an efficient reduction of the problem for general nonnegative matrices to that of doubly stochastic matrices. This was done via matrix (Sinkhorn's) scaling: for any matrix $A = \(a_{ij}\)$ with nonnegative entries and positive permanent, one can efficiently find {\it scaling factors} $x_1,\ldots,x_n$ and $y_1,\ldots,y_n$ such that the matrix $B = \(x_i \cdot a_{ij} \cdot y_j\)$ is (almost) doubly stochastic. Since $Per(A) = \frac{1}{\prod_i x_i \cdot \prod_j y_j} \cdot Per(B)$ this constitutes a reduction, and in fact achieves $e^n$ deterministic approximation for the permanent of a nonnegative matrix.

\subsection{Our results}
Our paper is a contribution to the third approach. One may say that, in a sense, it takes up where \cite{LSW} left off. The algorithm of \cite{LSW} reduces the problem to the case of doubly stochastic matrices, on which it "does nothing", that is returns $1$ and quits. The natural next step would be to "actually look at the matrix", that is to come up with an efficiently computable function of the entries of the matrix, which would provide a non-trivial estimate of its permanent.

\noi This is precisely what we do. This efficiently computable function of the doubly stochastic matrix $A = \(a_{ij}\)$ is
$F(A) = \prod_{i,j = 1^n} \(1-a_{ij}\)^{1-a_{ij}}$.

\noi We prove new lower and upper bounds for the permanent of a doubly stochastic matrix $A$, showing that for any such matrix it holds that
\beqn
\label{perm:bounds}
F(A) \le Per(A) \le 2^n \cdot F(A)
\eeqn
Combined with the preceding discussion, this gives our main algorithmic result.
\thm
\label{thm:approx}
There is a deterministic polynomial-time algorithm to approximate the permanent of a nonnegative matrix up to a multiplicative factor of $2^n$.
\ethm

\noi Let us now briefly describe the ideas leading to the bounds in (\ref{perm:bounds}).

\noi We proceed via convex relaxation. That is, given a matrix $A$ with nonnegative entries, we define a concave maximization problem, whose solution approximates $\log(Per(A))$.

\noi Let us start with pointing out that approximating the permanent via matrix scaling may
also be achieved by solving a convex optimization problem. In fact, what we need is to find the product of scaling factors $\prod_i x_i \cdot \prod_j y_j$ of $A$. This could be done in two different ways:

\noi By solving a concave maximization problem:
\beqn
\label{kld}
\log \(\frac{1}{\prod_i x_i \cdot \prod_j y_j}\) = \max_{B \in \Omega_n} \sum_{1 \leq i,j \leq i,j} b_{i,j} \log\(\frac{ a_{i,j}}{b_{i,j}}\)
\eeqn
Here $\Omega_n$ is the set of all $n\times n$ doubly stochastic matrices.

\noi And by solving a convex minimization problem:
\beqn \label{prd-pol}
\log \(\frac{1}{\prod_i x_i \cdot \prod_j y_j}\) = \inf_{x_1+...+x_n = 0} \log\(Prod_{A}\(e^{x_1},...,e^{x_n}\)\),
\eeqn
where $Prod_{A}(x_1,...,x_n)$ is the {\it product polynomial} of $A$,
\[
Prod_{A}\(x_1,...,x_n\) = \prod_{1 \leq i \leq n} \sum_{1 \leq j \leq n} a_{ij} x_j
\]

\noi Note that $Per(A)$
is the {\it mixed derivative} of $Prod_{A}$: $Per(A)= \frac{\partial^{n}}{\partial x_{1}\dots \partial x_n}Prod_{A}(0,...,0)$.

\noi The relaxation (\ref{kld}) is very specifically tied to the permanent. On the other hand, (\ref{prd-pol}) is much more general, in that it aims to approximate the mixed derivative of a homogeneous polynomial $p(x_1,...,x_n)$ of degree $n$ with non-negative coefficients, given via an evaluation oracle\footnote{Note that the product polynomial can be efficiently evaluated.}.

\noi In \cite{Gurvits-S}, the relaxation (\ref{prd-pol}) was shown to provide an $e^n$-approximation of the mixed derivative for a large class of homogeneous polynomials, containing the product polynomial. Moreover, it is the first step in a hierarchy of sharper relaxations given by considering
\[
\gamma_i =: \inf_{x_1+...+x_i = 0} \log(Q_i(e^{x_1},...,e^{x_n})),
\]
where $Q_i(x_1,...,x_i) = \frac{\partial^{n-i}}{\partial x_{i+1}\dots \partial x_n}p(x_1,...,x_i,0,...,0)$.

\noi If the (multivariate) polynomial $p$ does not have roots with positive real parts (in this case it is known as {\it H-Stable}, or {\it hyperbolic}) then
$$
G\(deg_{Q_{i+1}}(i+1)\) \cdot exp(\gamma_i) \leq Per(A) \leq exp(\gamma_i) \leq exp(\gamma_{i+1}),
$$
where $G(k) =: \left(\frac{k-1}{k}\right)^{k-1}$ and $deg_{Q_{i+1}}(i+1)$ is the degree of the variable $x_{i+1}$ in the polynomial
$Q_{i+1}$. In particular,
\beqn \label{her}
\frac{i!}{i^i} \cdot exp(\gamma_i) \leq Per(A) \leq exp(\gamma_i).
\eeqn
Considering this hierarchy turns out to be very useful, both from mathematical and from algorithmic points of view \cite{Gurvits-S}, \cite{L-S}. Note that, when this approach is applied to the product polynomial $Prod_{A}$, the original matrix structure is essentially lost. But by giving up the matrix structure, we gain additional
inductive abilities, leading, in particular, to a rather simple proof of (\ref{her}).

\noi Unfortunately, we only know how to compute $\gamma_i$ in $poly(n) \cdot 2^{n-i}$ oracle calls, which is polynomial-time only for $i = n - O(\log(n))$. In other words, this "hyperbolic polynomials" approach
does not seem to break the $e^n$-barrier for the approximation of the permanent by a polynomial-time deterministic algorithm.

\noi So, the challenge was to come up with a better
convex relaxation. Such a relaxation was suggested in \cite{Chertkov}, and it is a generalization of (\ref{kld}). It is a special case of a well-known heuristics in Machine Learning, the so called {\it Bethe Approximation}. This heuristics is used to approximate {\it log partition functions} of the following type (appearing, in particular, in the analysis of {\it Belief Propagation} algorithms).
\beqn \label{graph-mod}
PF =: \log\(\sum_{\begin{array}{l} x_i \in S_i \\ i = 1...n\end{array}} \prod_{i} G_i\(x_i\) \cdot \prod_{(i,j) \in E} F_{i,j}\(x_i,x_j\)\)
\eeqn
Here $S_i$ are finite sets; $G_i(x_i)$ and $F_{i,j}(x_i,x_j)$ are given non-negative functions, and $E$ is the set of edges of the associated undirected graph $\Gamma$.

\noi If the graph $\Gamma$ is a tree then $PF$ can be efficiently evaluated, e.g. by dynamic programming. The Bethe Approximation is a heuristic
to handle possible cycles. It turns out that $\log(Per(A)$ can be represented as in (\ref{graph-mod}). This was first observed in \cite{Jebara}. In this paper we use a simplified version of this heuristic proposed in \cite{Chertkov}, which amounts to approximating the logarithm of the permanent of a nonnegative matrix $A$ by
\beqn \label{ch-je}
\max_{B \in \Omega_n} \sum_{i,j=1}^n \(1- b_{ij}\) \log\(1- b_{ij}\)  + \sum_{i,j=1}^n b_{ij} \log\left(\frac{a_{ij}}{b_{ij}} \right).
\eeqn

\noi We should mention that, according to \cite{Heil}, the physicists had already applied the Bethe Approximation to the closely related monomer-dimer problem as early as in late 1930s.

\subsubsection*{Lower bound}
We prove that (\ref{ch-je}) is a lower bound on $\log(Per(A))$.
\thm
\label{thm:lower:general}
Let $A = \(a_{ij}\)_{i,j=1}^n$ be a nonnegative matrix and let $B = \(b_{ij}\)_{i,j=1}^n$ be a doubly stochastic matrix. Then
\beqn
\label{bnd:lower:general}
Per(A) \ge \prod_{i,j = 1}^n \(1-b_{ij}\)^{1-b_{ij}} \cdot \exp\left\{-\sum_{i,j=1}^n b_{ij} \log \frac{b_{ij}}{a_{ij}}\right\}
\eeqn
\ethm
Let us note that this claim was first stated (but not proved) in \cite{Vontobel}.

\noi If $A$ is doubly stochastic, setting $B = A$ in (\ref{bnd:lower:general}) gives the lower bound in (\ref{perm:bounds}).

\noi Theorem~\ref{thm:lower:general} has an additional combinatorial application. We show it to imply S. Friedland's "Asymptotic Lower Matching Conjecture" for the monomer-dimer problem. We will go into details in Section~\ref{sec:lower}.

\subsubsection*{Upper bound}
We prove that $2^n$ times (\ref{ch-je}) is an upper bound on $\log(Per(A))$.
\thm
\label{thm:upper bound}
The permanent of a stochastic matrix $A = \(a_{ij}\)$ satisfies
\[
Per(A) \le C^n \cdot \prod_{ij} \(1-a_{ij}\)^{1-a_{ij}}
\]
for some $C \le 2$.
\ethm
Note that this implies, in particular, that for a nonnegative matrix $A$, and its {\it doubly stochastic scaling} $B$, we have
\[
Per(A) \le 2^n \cdot \prod_{i,j = 1}^n \(1-b_{ij}\)^{1-b_{ij}} \cdot \exp\left\{-\sum_{i,j=1}^n b_{ij} \log \frac{b_{ij}}{a_{ij}}\right\}
\]

\rem
\begin{itemize}
\item
Let
\[
CW(A,B) = \sum_{i,j=1}^n \(1- b_{ij}\) \log\(1- b_{ij}\)  + \sum_{i,j=1}^n b_{ij} \log\left(\frac{a_{ij}}{b_{ij}} \right)
\]
The functional $CW(A,B)$ is clearly concave in $A$. Less obviously, it is concave in $B \in \Omega_{n}$ \cite{Vontobel}. So, in principle, the concave maximization problem (\ref{ch-je}) can be solved in polynomial deterministic time by, say, the ellipsoid method.

\noi We don't use the concavity in $B$ in this paper. The algorithm we propose and analyze first scales the matrix $A$ to a doubly-stochastic matrix $D$ and outputs  $\prod_{i,j = 1}^n \(1-d_{ij}\)^{1-d_{ij}}$ multiplied by the product of the scaling factors. So, when applied to a doubly-stochastic matrix, our algorithm has linear complexity.

\noi There are several benefits in using this suboptimal algorithm. First: We can analyze it.
Second: It is fast, and local (looking only at the entries) in the doubly-stochastic case. Third: it already improves on $e^n$-approximation. Fourth: it might allow (conjectural) generalizations to the
hyperbolic polynomials setting, to be described in the journal version.

\noi We also conjecture that our algorithm, might in fact turn out to be optimal. That is, that its worst case accuracy is the same as that of the Bethe Approximation (\ref{ch-je}).

\item
Let us remark that our results can be viewed as reasonably sharp bounds on a specific partition function in terms of its Bethe Approximation. To the best of our knowledge, this might be one of the first results of this type, and one of the first applications of the Bethe Approximation to theoretical computer science.
\end{itemize}
\erem

\noi {\bf Discussion}. It would seem that the improvement of the approximation factor from one exponential to a smaller one leaves something to be desired. This is, of course, true. On the other hand, let us remark that any algorithm which considers only the distribution of the entries of the matrix cannot achieve better than $2^{n/2}$ approximation for the permanent. This was pointed out to us by \cite{Wigderson-personal}. In fact, consider the following two $0,1$ matrices, both having $2$ ones in each row and column. The matrix $A_1$ is a block-diagonal matrix, with $n/2$ blocks of $\small{\left(\begin{array}{ll}
1 & 1 \\
1 & 1
\end{array}\right)}$ on the diagonal (assume $n$ is even). The matrix $A_2$ is the adjacency matrix of a $2n$-cycle, viewed as a bipartite graph with $n$ vertices on each side. The permanent of $A_1$ is clearly $2^{n/2}$, while the permanent of $A_2$ is $2$.

\noi We conjecture that this optimal approximation factor of $2^{n/2}$ can be attained, by improving our upper bound.
\cnj
\label{cnj:upper}
The permanent of a doubly stochastic matrix $A = \(a_{ij}\)$ satisfies
\[
Per(A) \le 2^{n/2} \cdot \prod_{ij} \(1-a_{ij}\)^{1-a_{ij}}
\]
\vspace*{-2ex}\ecnj
Note that this conjectured bound would be tight for the doubly stochastic matrix $\frac12 \cdot A_1$.

\noi {\bf Organization}: The organization of this paper is as follows: We discuss known combinatorial bounds for the permanent and their relation to our bounds in Section~\ref{sec:bounds}. We prove the lower bound in Section~\ref{sec:lower}, and the upper bound in Sections~\ref{sec:upper}~and~\ref{app:sec:upper}.

\section{Bounds for the permanent}
\label{sec:bounds}
\subsection{Lower bounds}
In general, the permanent of a nonnegative matrix may vanish. Hence, we need to impose additional constraints on the matrix to allow non-trivial lower bounds. Usually, the matrix is assumed to be {\it doubly stochastic}, that is to have row and column sums equal $1$. In this case it is easy to see that the permanent has to be positive. The most famous bound for permanents is that of Egorychev \cite{Eg} and Falikman \cite{Fal}, resolving the question of van der Waerden, and showing the permanent of a doubly stochastic matrix to be at least $\frac{n!}{n^n}$. This bound is tight and is attained on the matrix all of whose entries equal $1/n$.

\noi If we impose additional constraints on the matrix, we may expect a stronger bound. The class $\Lambda(k,n)$ of integer matrices whose row and column sums equal $k$ (adjacency matrices of $k$-regular bipartite graphs with multiple edges) was considered by Schrijver and Valiant \cite{S-V}. Normalizing by $k$, one obtains a class of doubly stochastic matrices with entries of the form $\frac{m}{k}$ for integer $m$ (and hence, with support of size at most $k$ in each row and column). The authors conjectured the minimal permanent for this class to be at least $\((k-1)/k\)^{(k-1)n}$. This conjecture was proved in \cite{Schrijver}\footnote{Let us remark that the assumption on the rationality of the entries was removed in \cite{Gurvits-S}, making only the structure of the support matter.}. A more general bound from \cite{Schrijver}  will be of special interest to us: Let $B = \(b_{ij}\)$ be a doubly stochastic matrix, and let $A = \(b_{ij} \cdot \(1 - b_{ij}\)\)$. Then

\beqn
\label{bnd:Schr}
Per(A) \ge \prod_{i,j=1}^n \(1-b_{ij}\)
\eeqn
We observe, for future reference, that the matrix $B$ is replaced by a new matrix $A$, obtained by applying a concave function $\phi(t) = t(1-t)$ entry-wise to $A$. For this new matrix, an explicit, efficiently computable, lower bound on the permanent is given.

\noi All these bounds are very difficult technical results, some of them using advanced mathematical tools, such as the Alexandrov-Fenchel inequalities. Let us note that more general bounds (with easier proofs), implying all the results above, were given in \cite{Gurvits-S}, using the machinery of hyperbolic polynomials. The point we would like to make (for future comparison with the situation with upper bounds) is that the lower bounds for the permanent are hard to prove, but they are essentially optimal.

\noi We now consider a more general notion than the permanent. For an $n \times n$ matrix $A$, and $1 \le m \le n$, let
$Per_m(A)$ be the sum of permanents of all $m \times m$ submatrices of $A$. Note that if $A$ is a $0,1$ matrix, the permanent counts the perfect matchings of the corresponding bipartite graph, while $Per_m(A)$ counts all the matchings with $m$ edges. Friedland \cite{FRIEDLAND} stated a conjectured lower bound on $Per_m$ for the class $\Lambda(k,n)$ of integer matrices\footnote{This lower bound is complicated, we will state it explicitly below.}. This conjecture has significance in statistical physics and is a natural generalization of the Schrijver-Valiant conjecture. Partial results towards this conjecture were obtained in \cite{fried}.

\noi {\bf Our results:}

\noi We restate our lower bound Theorem~\ref{thm:lower:general} here for the convenience of the reader:

\noi {\it
Let $A = \(a_{ij}\)_{i,j=1}^n$ be a nonnegative matrix and let $B = \(b_{ij}\)_{i,j=1}^n$ be a doubly stochastic matrix. Then
\[
Per(A) \ge \prod_{i,j = 1}^n \(1-b_{ij}\)^{1-b_{ij}} \cdot \exp\left\{-\sum_{i,j=1}^n b_{ij} \log \frac{b_{ij}}{a_{ij}}\right\}
\]
}

\noi We note that this lower bound is the first lower bound on the permanent which actually "looks at the matrix", that is depends explicitly on the entries of $A$, rather than on its support pattern.

\noi Note that the bound (\ref{bnd:Schr}) follows, by taking $A = \(b_{ij} \cdot \(1 - b_{ij}\)\)$. Hence Theorem~\ref{thm:lower:general} is a generalization of (\ref{bnd:Schr}). On the other hand, let us say that we view it as a corollary of (\ref{bnd:Schr}), since it is proved by analysis of the first order optimality conditions on the RHS of the inequality above, viewed as a function on doubly stochastic matrices, and the key part of the analysis is applying (\ref{bnd:Schr}).

\noi {\bf The conjecture of Friedland}.
Let $\alpha(m,n,k) = \min_{A \in \Lambda(k,n)} Per_{m}(A)$. Think about $m$ growing linearly in $n$ and $k$ being fixed\footnote{The bounds below hold for any $k$, though.}. Then $\alpha(m,n,k)$ is exponential in $n$, and we are interested in the exponent.

\noi To be more precise, fix $p \in [0,1]$ (this is the so called {\it limit dimer density}).
Let $m(n) \leq n$ be an integer sequence with $\lim_{n \rightarrow \infty}\frac{m(n)}{n} = p$.
Finally, let\footnote{It follows from Theorem~\ref{thm:Friedland cnj} that this definition is independent of the choice of the sequence $m(n)$ and that the limit exists.}
\[
\beta(p,k) = \lim_{n \rightarrow \infty} \frac1n \log(\alpha(m(n),n,k))
\]

\noi The challenge is to find $\beta(p,k)$. S. Friedland had conjectured that, similarly to \cite{Schrijver}, one can replace the minimum in the definition of $\alpha(m,n,k)$ by an (explicitly computable) average over a natural distribution $\mu = \mu_{k,n}$ on $\Lambda(k,n)$ (see Section~\ref{sec:lower}).

\noi We show this conjecture to hold, deducing it from Theorem~\ref{thm:lower:general}.
\thm
\label{thm:Friedland cnj}
\[
\beta(p,k) = \lim_{n \rightarrow \infty} \frac 1n \log\(\E_{\mu}(Per_{m(n)}(A))\)
\]
\ethm

\rem
Friedland's conjecture was proved, using the hyperbolic polynomials, in \cite{fried} for limit dimer densities of the form $p = \frac{k}{k+s}, s \in \N$.
\erem

\subsection{Upper bounds}
The notable upper bound for the permanents is due to Bregman \cite{Bregman}, proving a conjecture of Minc. This is a bound for permanents of $0,1$ matrices. For a $0,1$ matrix $A$ with $r_i$ ones in the $i^{th}$ row,
\beqn
\label{bnd:Bregman}
Per(A) \le \prod_{i=1}^n \(r_i!\)^{1/r_i}
\eeqn

\noi To the best of our knowledge, there is no satisfying extension of this bound to general nonnegative matrices. We will now give a different view of (\ref{bnd:Bregman}), suggesting a natural way to extend it. Let $A = \(a_{ij}\)$ be a stochastic matrix, whose values in the $i^{th}$ row are either $0$ or $1/r_i$. Let $B = \(b_{ij}\)$ be a matrix with $b_{ij} = 0$ if $a_{ij} = 0$ and $b_{ij} = \(1/r_i!\)^{1/r_i}$ if $a_{ij} = 1/r_i$. Then: $Per(B) \le 1$.

\noi There is a natural construction of a function on the interval $[0,1]$ taking $1/r$ to $\(1/r!\)^{1/r}$ for all integer $r$. This is the function $\phi_0(x) = \Gamma\(\frac{1+x}{x}\)^{-x}$.
\cnj (\cite{Sam})
\label{cnj:Bregman}
Let $A = \(a_{ij}\)$ be a stochastic matrix, and let $B = \(\phi_0\(a_{ij}\)\)$. Then $Per(B) \le 1$.
\ecnj
Unfortunately, we do not know how to prove this conjecture.

\noi There is, however, a way to view it as a special (difficult) case in a general family of upper bounds for the permanent. The function $\phi_0(x) = \Gamma\(\frac{1+x}{x}\)^{-x}$ is a {\it concave} \cite{Soules} increasing function taking $[0,1]$ onto $[0,1]$. We can ask for which concave functions $\phi$ of this form, Conjecture~\ref{cnj:Bregman} holds. Note the similarity of this point of view with that of the bound (\ref{bnd:Schr}). In both cases we apply a concave function entry-wise to the entries of a stochastic matrix and ask for an explicit efficiently computable upper (or lower) bound for the permanent of the obtained matrix.

\noi Let $\phi$ be concave increasing function taking $[0,1]$ onto $[0,1]$. The function $\psi = \phi^{-1}$ is convex increasing taking $[0,1]$ onto $[0,1]$. It defines an {\it Orlicz norm} (\cite{Zygmund}) $\|\cdot\|_{\psi}$ on $\R^n$ as follows: for $v = \(v_1,\ldots,v_n\) \in \R^n$
\[
\|v\|_{\psi} = s, \quad \mbox{where $s$ is such that } \sum_{i=1}^n \phi\(\frac{|v_i|}{s}\) = 1
\]
Note that this is a generalization of the more familiar $l_p$ norms. For $\psi(x) = x^p$, $\|\cdot\|_{\psi} = \|\cdot\|_p$.

\noi If $v$ is a stochastic vector, the vector $w = \(\phi\(v_1\),\ldots,\phi\(v_n\)\)$ has $\|w\|_{\psi} = 1$. Thus, the question we are asking is: For which Orlicz norms $\|\cdot\|_{\psi}$, a matrix $B$ whose rows are unit vectors in this norm has permanent at most $1$. Using homogeneity of the norm and multilinearity of the permanent, we obtain an appealing form of the general family of upper bounds to consider: We want any nonnegative matrix $B$ with rows $b_1,\ldots,b_n$ satisfy
\beqn
\label{bnd:upper:gen}
Per(B) \le \prod_{i=1}^n \|b_i\|_{\psi}
\eeqn
{\bf Our results:} We prove (\ref{bnd:upper:gen}) for a family of functions $\psi$. Theorem~\ref{thm:upper bound} follows as a corollary.

\noi We note, that in strong contrast to the lower bounds case, our bounds are far from being optimal, and, in particular, are far from proving Conjecture~\ref{cnj:upper} or Conjecture~\ref{cnj:Bregman}.

\section{Proofs of the lower bounds}
\label{sec:lower}

\subsection{Proof of Theorem~\ref{thm:lower:general}}

{\bf Notation}. We will denote by $\Omega_n$ the class of doubly stochastic $n\times n$ matrices. For a pair $P = \(p_{ij}\), Q = \(q_{ij}\)$ of non-negative matrices, we let
\[
CW(P,Q) = \sum_{i,j=1}^n \(1- q_{ij}\) \log\(1- q_{ij}\)  - \sum_{i,j=1}^n q_{ij} \log\left(\frac{q_{ij}}{p_{ij}} \right)
\]
Let $P$ be a non-negative $n \times n$ matrix with positive permanent (which we may assume, without loss of generality). We will prove the theorem by showing
\[
\log(Per(P)) \geq max_{Q \in \Omega_{n}} CW(P,Q)
\]
\noi Note that, by continuity, we may assume all the entries in $P$ to be strictly positive.
Then the functional $CW(P,Q)$ is bounded from above and continuous as function of $Q$ on $\Omega_{n}$. Therefore, the maximum is attained. Let $V \in \Omega_{n}$ be one of points at which it is attained.

\noi We first isolate ones in the doubly-stochastic matrix $V$: up to rearrangement of the rows and columns, $V = \left( \begin{array}{cc}
		  I & 0 \\
		  0 & T \end{array} \right)$, where the doubly -stochastic matrix $T$ does not have ones;
and block-partition accordingly the matrix $P = \left( \begin{array}{cc}
		  P^{(1,1)} & P^{(1,2)} \\
		  P^{(2,1)} & P^{(2,2)} \end{array} \right)$.\\

\noi Note that $CW(P,V) = CW\(P^{(2,2)}, T\) + \sum_i \log\(P^{(1,1)}_{i,i}\)$.

\noi Since $Per(P) \geq Per\(P^{(1,1)}\) \cdot Per\(P^{(2,2)}\) \geq \prod_i P^{(1,1)}_{i,i} \cdot Per\(P^{(2,2)}\)$, we only need to prove $\log\(Per(P^{(2,2)})\) \geq CW\(P^{(2,2)}, T\)$.

\noi Let $d$ be the dimension of matrices $P^{(2,2)}, T$.
We express the local extremality conditions for $T$ not on the full $\Omega_d$ but rather in the interior of the
compact convex subset of doubly-stochastic $d \times d$ matrices supported on the support of $T = \(t_{kl}\)$.

\noi We first compute the partial derivatives (writing them out for general $d$-dimensional $P,Q$)
\[
\frac{\partial}{\partial q_{ij}} CW(P, Q) =  -2 - \log\(1- q_{ij}\) - \log\(q_{ij}\) + \log\(p_{ij}\) \quad 1 \leq i,j \leq d
\]
\noi By the first order optimality conditions for $T$, we get that there exists real numbers $\{\alpha_k\}, \{\beta_l\}$
such that
$$
-2 -\log\(1-t_{kl}\) - \log\(t_{kl}\) + \log\(P^{(2,2)}_{kl}\) = \alpha_k + \beta_l;~ (k,l) \in Supp(T)
$$
Which gives, for some positive numbers $\{a_k\}, \{b_l\}$ the following scaling:
\[
P^{(2,2)}_{kl} = a_k b_l \cdot t_{kl} \(1 - t_{kl}\);~ (k,l) \in Supp(T)
\]
Now, we can conclude the proof.
\begin{enumerate}
\item
It follows from the definition of the support that (applying the inequality below entry-wise)
\[
P^{(2,2)} \geq Diag\(a_k\) \cdot \widetilde{T} \cdot Diag(b_l); ~\mbox{where}~~\widetilde{T}_{kl} = t_{kl} \(1 - t_{kl}\)
\]
\item
It follows from doubly-stochasticity of $T$ that
\beqn
\label{dva}
CW(P^{(2,2)},T) = \sum \log\(a_k\) + \sum \log\(b_l\) + \sum_{(k,l) \in Supp(T)} \log\(1 - t_{kl}\)
\eeqn
\end{enumerate}
Finally it follows from (\ref{dva}) and (\ref{bnd:Schr}) that
$$
\log\(Per\(Diag\(a_k\) \cdot \widetilde{T} \cdot Diag\(b_l\)\)\) \geq CW\(P^{(2,2)},T\)
$$
and therefore
$$
\log\(Per\(P^{(2,2)}\)\) \geq \log\(Per\(Diag\(a_k\) \cdot \widetilde{T} \cdot Diag\(b_l\)\)\) \geq CW\(P^{(2,2)},T\)
$$
\eprf

\subsection{Proof of Theorem~\ref{thm:Friedland cnj}}
Let us first recall the following well known identity (see, for instance, \cite{friedland-tverberg}), expressing $Per_m(A)$ as a single permanent:
\[
Per_m(A) =((n-m)!)^{-2} \cdot Per (L), \quad L =  \left( \begin{array}{cc}
		  A & J_{n,n-m} \\
		  J_{n,n-m}^{T}&0 \end{array} \right)
\]
where $J_{n,n-m}$ is $n \times (n-m)$ matrix of all ones. If the matrix $A \in c \cdot \Omega_n$ (i.e. proportional to a doubly-stochastic matrix) then
it is easy to scale the matrix $L$. In particular, if $A \in \Lambda(k,n)$ then
\beqn \label{ba}
Per_{m}(A) = \frac{Per(K)}{a^m b^{2(n-m)} ((n-m)!)^2}
\eeqn
where $K \in \Omega_{2n - m}$ is defined as follows
\[
K =  \left( \begin{array}{cc}
		  a \cdot A & b \cdot J_{n,n-m} \\
		  (b \cdot J_{n,n-m})^{T}& 0 \end{array} \right)
\]	
with $p = \frac{m}{n}$, $a = \frac{p}{k} = \frac{m}{kn}$, $b = \frac{1}{n}$.

\noi We note that the identity (\ref{ba}) follows from the diagonal scaling:
\[
K =\left(\sqrt{a}I_n \oplus \frac{b}{\sqrt{a}}I_{n-m}\right) \cdot  L \cdot \left(\sqrt{a}I_n \oplus \frac{b}{\sqrt{a}}I_{n-m}\right)
\]

\noi To proceed with the proof, we will need the following simple claim, following from the convexity of $(1-x) \log(1-x)$.
\pro
\label{pro:friedland}
Let $p_1,...,p_k$ be non-negative numbers, with $0 \leq p_i \leq 1$ and $\sum_{i=1}^k p_i = s$. Then, setting $b = \frac{s}{k}$,
\[
\prod_{i=1}^k \(1- p_i\)^{1-p_i} \geq (1 - b)^{k(1-b)}
\]
\epro

\noi Our main claim is:
\thm
Let $A \in  \Lambda(k,n)$, Let $1 \leq m \leq n$ and let $p =\frac{m}{n}$ . Then the following inequality holds\footnote{Assuming, for typographic simplicity, all the relevant values on LHS to be integer.}:
\beqn \label{pa1}
Per_{m}(A) \geq \frac{(\frac{k-p}{k})^{n(k-p)} \cdot (1 - n^{-1})^{(1 - n^{-1})2n^2(1-p)}}{(\frac{p}{k})^{np} \cdot n^{-2n(1-p)} \cdot ((n(1-p))!)^2}
\eeqn
\ethm
 \prf
Apply the lower bound in (\ref{perm:bounds}) to the doubly-stochastic matrix $K$ and use (\ref{ba}). If $A$ is boolean then this already gives the inequality we need. In the non-boolean case an immediate application of Proposition~\ref{pro:friedland} finishes the proof.
\eprf

\noi {\bf Proof of Theorem~\ref{thm:Friedland cnj}}.

\noi First, we define the distribution $\mu$ on $\Lambda(k,n)$. Consider the following construction of a matrix $A \in \Lambda(k,n)$. For a permutation $\pi \in S_{kn}$, let $M = M_{\pi}$ be the standard representation of $\pi$ as a $kn \times kn$ matrix of zeroes and ones. Now, view $M$ in the natural way as a $k \times k$ block matrix $M = \(M_{ij}\)$, where each block $M_{ij}$ is an $n \times n$ matrix. Finally, set $A = A(\pi) = \sum_{i,j=1}^k M_{ij}$. The distribution $\mu$ is the one induced on $\Lambda(k,n)$ by the uniform distribution on $S_{kn}$.

\noi We point out that the expectation $\E_{\mu} \(Per_{m}(A)\)$ is known (see for instance \cite{FRIEDLAND}, \cite{fried}). In particular,
if $\lim_{n \rightarrow \infty} \frac{m(n)}{n} = p \in [0,1]$ then the following equality holds:
\beqn \label{ira}
\lim_{n \rightarrow \infty} \frac{\log\(\E_{\mu} \(Per_{m(n)}(A)\)\)}{n} = p \log\left(\frac{k}{p}\right) -2(1-p)\log(1-p) + (k-p) \log\(1 - \frac{p}{k}\)
\eeqn
The claim of the theorem follows directly from (\ref{pa1}), (\ref{ira}), and Stirling's formula.
\eprf

\section{Proofs of the upper bounds}
\label{sec:upper}
Recall that we are interested in upper bounds of the form given in (\ref{bnd:upper:gen}). We prove the following general claim.
\thm
\label{thm:general}
Let $\psi$ be a convex increasing thrice differentiable function taking $[0,1]$ onto $[0,1]$. Assume $\psi$ has the following properties
\begin{enumerate}
\item
The function $x \cdot \frac{\psi'(x)}{\psi(x)}$ is increasing.
\item
The function $x \cdot \frac{\psi''(x)}{\psi'(x)}$ is increasing.
\item
\[
\psi\(e^{-r/e}\) + \psi\(r \cdot e^{-r/e}\) \ge 1~~~~\mbox{for $0 \le r \le 1$}
\]
\end{enumerate}

\noi Then, for any nonnegative matrix $B$ with rows $b_1,\ldots,b_n$ it holds that
\[
Per(B) \le \prod_{i=1}^n \|b_i\|_{\psi}
\]
\ethm

\noi For this theorem to be useful, we need to provide examples of functions it applies to. We now give an example of a function $\psi$ satisfying the conditions of the theorem. Let $a \approx 1.54$ be the unique root of the equation
$
\frac{1-\ln a}{a} = \frac1e
$.

\lem
\label{lem:psi-exam-good}
The function
\[
\psi_a(x) = 1 - (1-x) \cdot a^x
\]
satisfies the conditions of Theorem~\ref{thm:general}.
\elem

\noi We now show how to deduce Theorem~\ref{thm:upper bound} from Theorem~\ref{thm:general}, using the function $\psi_a$. We start with a technical lemma.
\lem
\label{lem:technical}
\begin{itemize}
\item
For any stochastic vector $x = \(x_1,\ldots,x_n\)$, the maximum of the entries of the vector $\(\frac{x_j}{\prod_{k=1}^n \(1-x_k\)^{1-x_k}}\)_{j=1}^n$ is at most $e^{1/e} \approx 1.44$.
\item
Let $\psi_a$ be the function in Lemma~\ref{lem:psi-exam-good}. Then for any stochastic vector $x = \(x_1,\ldots,x_n\)$ holds\footnote{Note that by the first claim of the lemma, all the arguments of $\psi$ in LHS are in the allowed range $[0,1]$.}
\[
\sum_{j=1}^n \psi_a\(\frac{x_j}{C \cdot \prod_{k=1}^n \(1-x_k\)^{1-x_k}}\) \le 1
\]
for some $e^{1/e} \le C \le 2$.
\end{itemize}
\elem
Given the lemma, Theorem~\ref{thm:upper bound} follows immediately: In fact, by the definition of $\|\cdot\|_{\psi}$, we have for any stochastic vector $x$,
\[
\|x\|_{\psi_a} \le C \cdot \prod_{k=1}^n \(1-x_k\)^{1-x_k}
\]
Hence, by Theorem~\ref{thm:general}, for any stochastic matrix $B$, whose rows are stochastic vectors $b_1,\ldots,b_n$,
\[
Per(B) \le \prod_{i=1}^n \|b_i\|_{\psi_a} \le C^n \cdot \prod_{i,j=1}^n \(1-b_{ij}\)^{1-b_{ij}}
\]
giving Theorem~\ref{thm:upper bound}.

\noi The full proofs of the claims in this section are given in the next section.

\section{Full proofs of the claims for the upper bound}
\label{app:sec:upper}

\subsection{Proof of Theorem~\ref{thm:general}}
A word on notation. We denote by $\|x\|_{\psi}$ the norm of a vector $x$ in $\R^k$, without stating $k$ explicitly. Thus, we may and will compare $\|\cdot\|_{\psi}$-norms of vectors of different dimensions.

\noi We denote by $A_{ij}$ the submatrix of a matrix $A$ obtained by removing the $i^{th}$ row and the $j^{th}$ column of $A$.

\noi The proof is by induction on the dimension $n$. For $n=1$ the claim holds since for a scalar $a \in \R$,
\[
Per(a) = a = \|a\|_{\psi}
\]
The second equality is due to the fact that $\psi(1) = 1$.

\noi Assume the theorem holds for $n-1$. The induction step from $n-1$ to $n$ is incorporated in the following lemma.
\lem
\label{lem:induction step}
Let $\g:~\R_+ \rarrow \R_+$ be a scalar function defined by
\[
\g(r) = \min_{y \in \R^{n-1}_+:~\|y\|_{\psi} = 1} \|(y,r)\|_{\psi}
\]
Assume $\g$ satisfies the following functional inequality: For any $r_1,\ldots,r_n \in \R_+$
\beqn
\label{perm:func-ineq}
\prod_{k=1}^n \g\(r_k\) \ge \sum_{k=1}^n r_k
\eeqn
Then, if the theorem holds for $n-1$, it holds also for $n$.
\elem
\prf {\bf of Lemma~\ref{lem:induction step}}

\noi Write the rows of the $n \times n$ matrix $A$ as $a_k = \(x_k,b_k\)$, with $x_k \in \R^{n-1}$ and $b_k = a_{kn} \in \R$.

\noi Clearly, if any of $a_k$ is $0$ the claim of the theorem holds. The other boundary case we need to treat separately is the case in which one of the vectors $x_k$ is $0$. Without loss of generality, assume $x_1 = 0$. Expanding the permanent with respect to the first row, and using the induction hypothesis for $A_{1n}$, we have
\[
Per(A) = a_{1n} \cdot Per\(A_{1n}\) \le a_{1n} \cdot \prod_{k=2}^n \|x_k\|_{\psi} \le \prod_{k=1}^n \|a_k\|_{\psi}
\]
establishing the theorem in this case.

\noi Assume none of $x_k$ is $0$. Expanding the permanent of $A$ with respect to the last column, and using the induction hypothesis, we have
\[
Per(A) = \sum_{i=1}^n b_i \cdot Per\(A_{in}\) \le \sum_{i=1}^n b_i \cdot \prod_{j \not = i} \|x_j\|_{\psi} =
\prod_{j=1}^n \|x_j\|_{\psi} \cdot \sum_{i=1}^n \frac{b_i}{\|x_i\|_{\psi}}
\]
Hence, to prove the theorem for $A$, we need to show
\[
\sum_{i=1}^n \frac{b_i}{\|x_i\|_{\psi}} \le \prod_{k=1}^n \frac{\|\(x_k,b_k\)\|_{\psi}}{\|x_k\|_{\psi}}
\]
Let $r_k = b_k/\|x_k\|_{\psi}$, $y_k = x_k/\|x_k\|_{\psi}$. Then the inequality translates to
\[
\prod_{k=1}^n \|\(y_k,r_k\)\|_{\psi} \ge \sum_{i=1}^n r_i
\]
which follows from (\ref{perm:func-ineq}), since $\|y_k\|_{\psi} = 1$, and hence $\|\(y_k,r_k\)\|_{\psi} \ge \g\(r_k\)$.
\eprf

\noi It remains to prove (\ref{perm:func-ineq}).

\noi First, we observe that the function $\g$ has an explicit form.
\lem
\label{lem:g-explicit}
\[
\g(r) = \|(1,r)\|_{\psi}
\]
\elem
\prf ({\bf of Lemma~\ref{lem:g-explicit}})

\noi We may assume $r > 0$, otherwise the claim of the lemma holds trivially.

\noi Consider the optimization problem of minimizing $\|(y,r)\|_{\psi}$ for $y$ in the unit sphere of the norm in $\R^{n-1}$. Note that the minimum is attained, since we are looking for the minimum of a continuous function in a compact set.

\noi Let $\y$ be a point of minimum. We will show $\y$ to be a unit vector, implying the claim of the lemma.

\noi {\bf First step:} We show $\y$ to be constant on its support.

\noi Since $\|(\y,r)\|_{\psi} = \g(r)$, we have
\[
\Big |\Big|\(\frac{\y}{\g(r)},\frac{r}{\g(r)}\)\Big |\Big |_{\psi} = 1 \le \Big |\Big |\(\frac{y}{\g(r)},\frac{r}{\g(r)}\)\Big |\Big |_{\psi}
\]
for any $y$ of norm $1$. Therefore $\z = \frac{\y}{\g(r)}$ is a point of minimum of $\sum_{i=1}^{n-1} \psi\(z_i\)$ in the domain $D = \left\{z:~\|z\|_{\psi} = 1/\g(r)\right\}$.

\noi Consider this new optimization problem. Set $a = \g(r)$ for typographic convenience. Note $a > 1$, since, by assumption, $r > 0$.  Then
\[
D = \left\{z \in R^{n-1}_+:~\sum_{i=1}^{n-1} \psi\(az_i\) = 1\right\}
\]
We know that $\z$ is a point of minimum of the target function $\sum_{i=1}^{n-1} \psi\(z_i\)$ on $D$.

\noi Let $S = S\(\z\)$ be the support of $\z$. The first order optimality conditions for $\z$ imply that there exists a constant $\l \in \R$ such that for any $i \in S$,
\beqn
\label{ineq-phi-der}
\frac{\psi'\(z_i\)}{\psi'\(az_i\)} = \lambda \cdot a
\eeqn

\noi We would like to deduce from this that $\z$ (and hence also $\y$) is constant on its support $S$.

\noi Let $\eta(x) = \ln \psi'\(e^x\)$. We claim that $\eta$ is strictly convex on $(-\infty,0]$. In fact, $\eta'(x) = \frac{e^x \psi''\(e^x\)}{\psi'\(e^x\)}$, which is strictly increasing in $x$, by the second assumption of the theorem.

\noi Note that $\psi'(x) = \exp\left\{\eta(\ln x)\right\}$. Therefore (\ref{ineq-phi-der}) is equivalent to
\[
\eta\(\ln\(z_i\)\) - \eta\(\ln\(z_i\) + \ln(a)\) = \ln\(\lambda \cdot a\)
\]
And this can't hold for different values of $z_i$ if $\eta$ is strictly convex. This shows $\z$ is constant on $S$, completing the first step.

\noi {\bf Second step}: $|S| = 1$.

\noi Let $|S| = k$, for some $1 \le k \le n-1$.

\noi Since $\sum_{i \in S} \psi\(a \cdot \(\z\)_i\) = 1$ and $\z$ is constant on $S$, we have for all $i \in S$,\\ $\(\z\)_i~=~\(1/a\) \cdot \psi^{-1}\(1/k\)$. Therefore
\beqn
\label{ineq:z-const}
\sum_{i=1}^{n-1} \psi\(\(\z\)_i\) = k \cdot \psi\(\frac{\psi^{-1}\(1/k\)}{a}\)
\eeqn

\noi Consider the function $f(x) = \(1/x\) \cdot \psi\(\frac{\psi^{-1}\(x\)}{a}\)$. We will show this function to decrease on the interval $[0,1]$. This would imply the minimum over $k$ of LHS of (\ref{ineq:z-const}) is attained at $k = 1$, completing this step.

\noi Taking the first derivative, and denoting $\alpha = \psi^{-1}$, we need to verify for $x \in (0,1)$
\[
0 > f'(x) = -\frac{1}{x^2} \cdot \psi\(\frac{\alpha\(x\)}{a}\) + \frac{1}{x} \cdot \psi'\(\frac{\alpha\(x\)}{a}\) \cdot \frac{\alpha'(x)}{a}
\]
That is,
\[
\psi\(\frac{\alpha\(x\)}{a}\) > \frac xa \cdot \psi'\(\frac{\alpha\(x\)}{a}\) \cdot \alpha'(x)
\]
\[
\psi\(\frac{\alpha\(x\)}{a}\) \cdot \psi'(\alpha(x)) >  \frac xa \cdot \psi'\(\frac{\alpha\(x\)}{a}\)
\]
Since $x = \psi(\alpha(x))$, we want to show
\[
\frac{\psi'(\alpha(x))}{\psi(\alpha(x))} > \frac 1a \cdot \frac{\psi'\(\frac{\alpha\(x\)}{a}\)}{\psi\(\frac{\alpha\(x\)}{a}\)}
~~\Longleftrightarrow~~\alpha(x) \cdot \frac{\psi'(\alpha(x))}{\psi(\alpha(x))} > \frac{\alpha(x)}{a} \cdot \frac{\psi'\(\frac{\alpha\(x\)}{a}\)}{\psi\(\frac{\alpha\(x\)}{a}\)}
\]
That is, it suffices to show that $y \cdot \frac{\psi'(y)}{\psi(y)}$ increases in $y$, and this is true by the first assumption of the theorem.

\noi This completes the second step and the proof of Lemma~\ref{lem:g-explicit}.

\eprf

\noi As the next step towards the proof of (\ref{perm:func-ineq}), we give a sufficient condition for a function $g:~\R_+ \rarrow \R_+$ to satisfy the functional inequality stated in (\ref{perm:func-ineq}) for $\g$.
\lem
\label{perm:lem:suff-ineq}
If
\[
g(x) \ge \left\{\begin{array}{lll}
e^{x/e} & for & 0 \le x \le e\\
x       &            & otherwise
\end{array}\right.
\]
then $\prod_{k=1}^n g\(r_k\) \ge \sum_{k=1}^n r_k$.
\elem
\prf
Let $0 \le r_1 \le r_2 \le \ldots \le r_n$ be given, and assume $r_k < e$, $r_{k+1} \ge e$.

\noi First assume $k < n$. Write $y = \sum_{i=1}^k r_i$, $z = \sum_{j=k+1}^n r_j$. Clearly, $z \ge e$. Note that, by assumption,  \[
\prod_{j=k+1}^n g\(r _j\) \ge \prod_{j=k+1}^n r _j \ge \sum_{j=k+1}^n r_j = z
\]
We have
\[
\prod_{i=1}^n g\(r_i\) = \prod_{i=1}^k g\(r_i\) \cdot \prod_{j=k+1}^n g\(r_j\) \ge e^{1/e \cdot \sum_{i=1}^k r_i} \cdot z =  e^{y/e} \cdot z
\]
It remains to show $e^{y/e} \cdot z  \ge y + z$ for $z \ge e$. Since $e^x \ge x + 1$, we have
\[
e^{y/e} \ge y/e + 1 \ge \frac{y + z}{z}
\]
and we are done in this case.

\noi The other case to consider is $k = n$. Write $y = \sum_{i=1}^k r_i$. In this case we need to show $e^{y/e} \ge y$ for all $y \ge 0$. This again follows from the inequality $e^x \ge x + 1$, substituting $x = y/e - 1$.
\eprf

\noi To prove (\ref{perm:func-ineq}) and complete the proof of the theorem, it remains to verify
$\g(r) = \|(1,r)\|_{\psi}$ satisfies the assumptions of Lemma~\ref{perm:lem:suff-ineq}. First, clearly,
\[
\g(r) \ge \|r\|_{\psi} = r
\]
Next, $\g(r) \ge e^{r/e}$ iff
\beqn
\label{perm:g-suff-cond}
\psi\(e^{-r/e}\) + \psi\(r \cdot e^{-r/e}\) \ge 1
\eeqn
So we need to verify this for $0 \le r \le e$.

\noi We now claim that we may reduce the problem to a subinterval.
\lem
\label{lem:subinterval}
Let $\psi$ be an increasing differentiable convex function, taking $[0,1]$ to itself. If $\psi\(e^{-r/e}\) + \psi\(r \cdot e^{-r/e}\) \ge 1$ on $[0,1]$, then this also holds for $[0,e]$.
\elem

\noi Observe that the third assumption of the theorem is that (\ref{perm:g-suff-cond}) holds for $r \in [0,1]$. Thus, proving the lemma will complete the proof of the theorem.

\prf
Set
\[
h(r) = \psi\(e^{-r/e}\) + \psi\(r \cdot e^{-r/e}\)
\]
Then
\[
h'(r) = \(e^{-r/e} - \frac1e r e^{-r/e}\) \cdot \psi'\(r e^{-r/e}\) - \frac1e e^{-r/e} \cdot \psi'\(e^{-r/e}\)
\]
First, we claim that $h'$ is nonnegative on $[1,e-1]$. In fact, on this interval $r e^{-r/e} \ge e^{-r/e}$. Consequently, by convexity of $\psi$, $\psi'\(r e^{-r/e}\) \ge \psi'\(e^{-r/e}\)$. Hence
\[
h'(r) \ge \psi'\(e^{-r/e}\) \cdot e^{-r/e} \cdot \(1 - \frac{r+1}{e}\) \ge 0
\]

\noi Next, we claim that $h(e-r) \ge h(r)$ for $0 \le r \le 1$. We need to show that
\[
\psi\((e-r) \cdot e^{-(e-r)/e}\) + \psi\(e^{-(e-r)/e}\) \ge \psi\(e^{-r/e}\) + \psi\(r \cdot e^{-r/e}\)
\]
Let $a, b$ be the arguments on LHS, and $c, d$ on RHS. Note $a \ge b$ and $c \ge d$. Since $\psi$ is convex and increasing, it will suffice to show $a+b \ge c+d$ and $a \ge c$ (this would imply $(a,b)$ majorizes $(c,d)$).
\begin{itemize}
\item
We argue $a+b \ge c+d$. Let $f(x) = (x+1) e^{-x/e}$, and let $g(x) = f(e-x)$. We want to show $g(x) \ge f(x)$ for $0 \le x \le 1$. Note that $f$ is increasing on $[0,e-1]$ and decreasing on $[e-1,e]$, so both $f$ and $g$ are increasing on $[0,1]$. First, we argue $f' \ge g'$. In fact, we have
\[
f'(x) = \frac1e \cdot ((e-1) - x) \cdot e^{-x/e} \ge g'(x) = \frac1e \cdot (1-x) \cdot e^{-(e-x)/e}
\]
So, it would suffice to check $g(1) \ge f(1)$ which, after simplification, is the same as $e^{1/e} \ge 2^{1/2}$. And this is true.
\item
We argue $a \ge c$, that is $(e-r) \cdot e^{-(e-r)/e} \ge e^{-r/e}$ on $[0,1]$. Let $g(x)$ be the first function, and $f(x)$ the second. Note that $f(0) = g(0) = 1$. Hence, it suffices to prove $f' \le g'$. We have $f'(x) = -1/e \cdot e^{-x/e}$ and $g'(x) = -e^{-(e-x)/e} + \frac{e-x}{e} \cdot e^{-(e-x)/e}$. Therefore
\[
g'(x) - f'(x) = \frac1e \cdot \((e-x) \cdot e^{-(e-x)/e} + e^{-x/e} - e\cdot e^{-(e-x)/e}\) =
\]
\[
\frac1e \cdot \(e^{-x/e} - x\cdot e^{-(e-x)/e}\) \ge 0
\]

\end{itemize}

\eprf

\subsection{Proof of Lemma~\ref{lem:psi-exam-good}}
We will prove the lemma in greater generality, that is for all functions $\psi = \psi_a$, with $\frac1e \le \frac{1-\ln a}{a} < 1$.

\noi First, we compute the first three derivatives of $\psi$.
\[
\psi'(x) = (1-(1-x) \cdot \ln a) \cdot a^x
\]
\[
\psi''(x) = \ln a \cdot (2-(1-x) \cdot \ln a) \cdot  a^x
\]
\[
\psi'''(x) = \ln^2 a \cdot (3-(1-x) \cdot \ln a) \cdot  a^x
\]

We now prove the required properties of $\psi$.
\begin{enumerate}
\item
For $1 < a < e$, the function $\psi$ is increasing strictly convex taking $[0,1]$ to $[0,1]$. In fact, by observation, $\psi' > 0$ for $0 \le x \le 1$ and $\psi'' > 0$ for $0 \le x \le 1$.

\item
The function $x \cdot \frac{\psi'(x)}{\psi(x)}$ is strictly increasing for $1 < a < \sqrt{e}$.\footnote{It is easy to check that all $a$ for which $\frac 1e \le \frac{1-\ln a}{a} <1$ lie in this interval.}

\noi It suffices to show for $0 < x < 1$
\[
\(\psi' + x\psi''\) \cdot \psi > x\(\psi'\)^2
\]
For typographic convenience, write $b = \ln a$. Substituting the expressions for $\psi$ and its derivatives, and introducing notation
\[
P(x) = b^2 x^2 + \(2b - 2b^2\)x + (1-b)^2,~~~~~Q(x) = b^2 x^2 + \(3b - b^2\)x + (1-b),
\]
we need to verify
\[
Q(x) \cdot \(1-(1-x)\cdot e^{bx}\) > xP(x)\cdot e^{bx}
\]
Observe that $Q$ is strictly positive on $(0,1)$. Rearranging, we need to show
\[
e^{-bx} > x \cdot \frac{P(x)}{Q(x)} + (1-x) = 1 - x \cdot \frac{Q(x) - P(x)}{Q(x)}
\]
Since $e^{-bx} > 1- bx$ on $(0,1)$, it suffices to show $(Q-P)/Q \ge b$, that is $(1-b)\cdot Q \ge P$. And this is directly verifiable, for $x \in (0,1)$ and $b \in \(0,1/2\)$.

\item
The function $x \cdot \frac{\psi''(x)}{\psi'(x)}$ is strictly increasing for $1 < a < \sqrt{e}$.

\noi This is true iff
\[
\(\psi''(x) + x \psi'''(x)\) \cdot \psi'(x) > x \cdot \(\psi''(x)\)^2
\]
Since $\psi''' > 0$, it suffices to prove
\[
\psi''(x) \cdot \psi'(x) \ge x \cdot \(\psi''(x)\)^2~~~\Longleftrightarrow~~~x \cdot \psi''(x) \le \psi'(x)
\]
Substituting the expressions for the derivatives of $\psi$ and simplifying, we need to verify
\[
bx(2-(1-x)b) \le 1 - (1-x)b
\]
This is a quadratic inequality in $x$. For $0 < b < 1/2$, the interval between the roots of this quadratic is easily seen to contain $[0,1]$, and we are done.

\item
\[
\psi\(e^{-r/e}\) + \psi\(r \cdot e^{-r/e}\) \ge 1~~~~\mbox{for $0 \le r \le 1$}
\]
As in the proof of Lemma~\ref{lem:subinterval}, we set
\[
h(r) = \psi\(e^{-r/e}\) + \psi\(r \cdot e^{-r/e}\)
\]
Hence
\[
h'(r) = \(e^{-r/e} - \frac1e r e^{-r/e}\) \cdot \psi'\(r e^{-r/e}\) - \frac1e e^{-r/e} \cdot \psi'\(e^{-r/e}\)
\]
Observe $h(0) = 1$. Hence, it suffices to prove $h' \ge 0$ on $[0,1]$. Equivalently, for $0 \le r \le 1$,
\[
\frac{\psi'\(r\cdot e^{-r/e}\)}{\psi'\(e^{-r/e}\)} \ge \frac{1}{e-r}
\]
Set $y = e^{-r/e}$. Clearly $e^{-1/e} \le y \le 1$. We will show a stronger statement
\[
\frac{\psi'\(ry\)}{\psi'\(y\)} \ge \frac{1}{e-r}
\]
for all $y$ in the range. Similarly to the argument in the first step in the proof of Lemma~\ref{lem:g-explicit}, $\ln\(\psi'\(e^x\)\)$ is convex in $x$, which implies the LHS is decreasing in $y$, so it suffices to prove the inequality for $y = 1$. Substituting the expression for $\psi'$ and again writing $b$ for $\ln a$, we need to verify
\[
(e-r) \cdot (1 - (1-r)b) \ge e^{b(1-r)},
\]
for $0 \le r \le 1$. At $r = 0$, we need to check $e \ge e^b / (1-b) = a/(1-\ln a)$, which is satisfied with equality, by the assumption. Clearly, RHS decreases in $r$. By a direct calculation, the derivative of LHS is positive, that is LHS is increasing, completing the proof.

\end{enumerate}

\subsection{Proof of Lemma~\ref{lem:technical}}

\noi For the first claim, we need a technical lemma.
\lem
\label{lem:gurvits-cnj:tech}
Let $x = \(x_1,\ldots,x_n\)$ be a stochastic vector. Let $y = x_1$. Then
\[
\prod_{k=1}^n \(1-x_k\)^{1-x_k} \ge \frac{(1-y)^{1-y}}{e^{1-y}}
\]
\elem
\prf
We need to show
\[
\prod_{k=2}^n \(1-x_k\)^{1-x_k} \ge e^{y-1} \quad \Longleftrightarrow \quad \sum_{k=2}^n \(1-x_k\) \ln\(1-x_k\) \ge y-1
\]
for nonnegative $x_2,\ldots,x_n$ summing to $a := 1-y$. Let $\x$ be minimizer of $f\(x_2,\ldots,x_n\) = \sum_{k=2}^n \(1-x_k\) \ln \(1-x_k\)$ on this domain. Let $S$ be the support of $\x$. The first order regularity conditions state the existence of a constant $\l$ such that
\[
\ln\(1-(\x)_k\) = \l
\]
for all $k \in S$. This means that $(\x)_k$ are constant on $S$.

\noi Let $s = |S|$. Then $f\(\x\) = (s-a) \ln\(\frac{s-a}{s}\)$. It remains to argue
\[
(s-a) \ln\(\frac{s-a}{s}\) \ge -a,
\]
for all integer $s \ge 1$. In fact, the function $g(s) = (s-a) \ln\(\frac{s-a}{s}\)$ of the real variable $s$ is non-increasing on $[1,\infty)$, since
$g'(s) = \ln\(1-a/s\) + a/s \le 0$. And it is easy to see that $g(s)$ tends to $-a$ as $s \rarrow \infty$.

\eprf

\noi This means that
\[
\frac{x_1}{\prod_{k=1}^n \(1-x_k\)^{1-x_k}} \le \frac{y e^{1-y}}{(1-y)^{1-y}}
\]
The following lemma concludes the proof of the first claim of Lemma~\ref{lem:technical}.
\lem
\label{lem:upper bound for f}
The function $f(y) = \frac{y e^{1-y}}{(1-y)^{1-y}}$ on $[0,1]$ is upperbounded by $e^{1/e}$.
\elem
\prf
The maximum $y e^{1-y}$ on $[0,1]$ is $1$ and the minimum of $(1-y)^{1-y}$ on $[0,1]$ is $e^{-1/e}$.

\eprf

\noi We move to the second claim of Lemma~\ref{lem:technical}, repeating its claim for convenience. Let $\psi$ be the function in Lemma~\ref{lem:psi-exam-good}. Then for any stochastic vector $x = \(x_1,\ldots,x_n\)$ holds
\[
\sum_{j=1}^n \psi\(\frac{x_j}{2 \cdot \prod_{k=1}^n \(1-x_k\)^{1-x_k}}\) \le 1
\]

\noi The proof contains two steps, given in the following lemmas.
\lem
\label{lem:gurvits-cnj:s1}
Let a stochastic vector $x = \(x_1,\ldots,x_n\)$ be given, and let $y = \max_i x_i$ be its maximal coordinate. Then, for any convex increasing function $\psi$ taking $[0,1]$ to itself, and for any constant $C \ge e^{1/e}$ it holds that
\beqn
\label{ineq:gurvits-cnj}
\sum_{j=1}^n \psi\(\frac{x_j}{C \cdot \prod_{k=1}^n \(1-x_k\)^{1-x_k}}\) \le \frac1y \cdot \psi\(\frac{y e^{1-y}}{C \cdot (1-y)^{1-y}}\)
\eeqn
\elem

\lem
\label{lem:gurvits-cnj:s2}
Let $\psi$ be the function in Lemma~\ref{lem:psi-exam-good}. Then
\[
\frac1y \cdot \psi\(\frac{y e^{1-y}}{2 \cdot (1-y)^{1-y}}\) \le 1
\]
for $0 < y \le 1$.
\elem
It remains to prove the lemmas.

\noi {\bf Proof of Lemma~\ref{lem:gurvits-cnj:s1}}

\noi Let $x = \(x_1,\ldots,x_n\)$ be a stochastic vector with maximal entry $y$. By Lemmas~\ref{lem:gurvits-cnj:tech}~and~\ref{lem:upper bound for f}, all the arguments of $\psi$ in the LHS of (\ref{ineq:gurvits-cnj}) are upperbounded by $m = \frac{y e^{1-y}}{C \cdot (1-y)^{1-y}} \le 1$ and their sum is at most $\frac{e^{1-y}}{C \cdot (1-y)^{1-y}}$. Since $\psi$ is convex and increasing, the maximum of LHS under these constraints is attained when $s = \lfloor 1/y \rfloor$ of these arguments equal $m$ and the remaining non-zero one equals $\frac{e^{1-y}}{C \cdot (1-y)^{1-y}} - s\cdot m$, which gives
\beqn
\label{a7:first}
\sum_{j=1}^n \psi\(\frac{x_j}{C \cdot \prod_{k=1}^n \(1-x_k\)^{1-x_k}}\) \le s \cdot \psi\(\frac{y e^{1-y}}{C \cdot (1-y)^{1-y}}\) + \psi\(\frac{(1 - sy) e^{1-y}}{C \cdot (1-y)^{1-y}}\)
\eeqn
Let $t = (1-sy)/y$. Since $s = \lfloor 1/y \rfloor$, we have $t \le 1$. Since $\psi$ is convex, increasing, and $\psi(0) = 0$, we have that $\psi(t \cdot x) \le t \cdot \psi(x)$ for any $t \le 1$, $0 \le x \le 1$. Therefore
\[
\psi\(\frac{(1 - sy) e^{1-y}}{C \cdot (1-y)^{1-y}}\) \le \frac{1-sy}{y} \cdot \psi\(\frac{y e^{1-y}}{C \cdot (1-y)^{1-y}}\)
\]
and the RHS of (\ref{a7:first}) is at most
\[
\frac1y \cdot \psi\(\frac{y e^{1-y}}{C \cdot (1-y)^{1-y}}\),
\]
completing the proof of the lemma.
\eprf

\noi {\bf Proof of Lemma~\ref{lem:gurvits-cnj:s1}}

\noi We have $\psi(x) = 1 - (1-x)\cdot a^x$, where $a \approx 1.54$ is determined by the identity $(1-\ln a)/a = \frac1e$.

\noi Set $f(y) = \frac1y \cdot \psi\(\frac{y e^{1-y}}{2 \cdot (1-y)^{1-y}}\)$.

\noi First, we claim that the maximum of $f$ on $[0,1]$ is attained for $y \le y_0 = 0.51$.

\noi In fact, setting $A(y) = \frac{y e^{1-y}}{C \cdot (1-y)^{1-y}}$, we have $f(y) = \psi(A(y))/y$ and $f'(y) \le 0$ iff
\[
\frac{A(y) \cdot \psi'(A(y))}{\psi(A(y))} \le \frac{1}{1 + y \ln(1-y)}
\]
Recall that the function $(x\cdot \psi'(x))/\psi(x)$ is increasing on $[0,1]$. Therefore the maximum of LHS is at most $a = \psi'(1)$. We claim that RHS is at least that, for $0.51 \le y \le 1$. In fact, this is easy to see that RHS is an increasing function of $y$. Computing this function at $y_0 = 0.51$, we see that it is greater than $a$, and we are done.
Hence $f'(0) < 0$ for $y_0 \le y \le 1$, and the maximum of $f$ is attained outside this interval.

\noi It remains to show $f(y) \le 1$ on $I = [0,y_0]$. Equivalently,
\beqn
\label{ineq-tech-gurvits}
1 - (1-A(y))\cdot a^{A(y)} \le y \quad \Longleftrightarrow \quad \ln(a) \cdot A(y) + \ln(1 - A(y)) \ge  \ln(1-y)
\eeqn
The function $\ln(a) \cdot x + \ln(1-x)$ is decreasing in $x$, and therefore we decrease LHS by substituting a larger value for $A(y)$. We prove (\ref{ineq-tech-gurvits}) for $y \in I$ by covering $I$ with several intervals, and, in each interval, replacing $A(y)$ by a different linear function which majorizes it in this interval.

\noi First, we need a technical lemma.
\lem
The function $g(y) = (1/y) \cdot \ln\(\frac{1-ry}{1-y}\)$ decreases on $(0,1/r)$ for any $r > 1$.
\elem
\prf
We will show $g' \le 0$. Computing the derivative and simplifying, we need to show
\[
(1-ry) \cdot \ln\(\frac{1-y}{1-ry}\) \le (r-1) \cdot \frac{y}{1-y}
\]
Since $\ln(1+x) \le x$, we may replace the logarithm on LHS with $\frac{(r-1)y}{1-ry}$, leading to a trivially true inequality.
\eprf

\noi We now prove (\ref{ineq-tech-gurvits}) in several steps. Observe, for future use, that the function $h(y) = e^{1-y}/(1-y)^{1-y}$ decreases on $[0,1]$.
\begin{itemize}
\item
The maximum of $h(y)$ on $[0,1]$ is $e = h(0)$. Set $r = e/2$. Then $A(y) \le ry$ for $y \in [0,1]$. (Note $ry < 1$ for $y \in I$.) Hence, if we show
\[
\ln(a) \cdot ry +  \ln\(1 - ry\) \ge \ln(1-y)
\]
for $y$ in some interval $I_1$, it would imply (\ref{ineq-tech-gurvits}) in this interval. Rearranging, we need to show
\[
\frac1y \cdot \ln\(\frac{1-ry}{1-y}\) \ge \ln(a) \cdot r
\]
By the lemma, LHS is a decreasing function of $y$, hence it suffices to check this inequality at the right endpoint of $I_1$. It holds at $y=0.3$, and therefore we may take $I_1 = [0,0.3]$ and (\ref{ineq-tech-gurvits}) holds in this interval.
\item
It remains to check (\ref{ineq-tech-gurvits}) in $[0.3,0.51]$. In this interval, the maximum of $h$ equals $h(0.3)$ and therefore $A(y) \le ry$ for $r = h(0.3)/2$. Repeating the same argument, with the new value of $r$, we extend the validity of (\ref{ineq-tech-gurvits}) to $I_2 = [0,0.4]$. Reiterating, with new values of $r$, we get, in two more steps, to progressively larger intervals $[0,0.48]$, and, finally, to $[0,0.51]$.
\end{itemize}
The lemma is proved.
\eprf


\begin{thebibliography}{99}
\bibitem{Barvinok}
A. I. Barvinok, {\sl Polynomial Time Algorithms to Approximate Permanents and Mixed Discriminants Within a Simply Exponential Factor.} Random Struct. Algorithms 14(1): 29-61 (1999)

\bibitem{Bregman}
L. M. Bregman, {\sl Certain properties of nonnegative matrices and their permanents},
Soviet Math. Dokl. 14, 945-949, 1973.

\bibitem{Broder}
A. Z. Broder, {\sl How hard is it to marry at random? (On the approximation of the permanent)}, in
Proceedings of the 18th Annual ACM Symposium on Theory of Computing (STOC), ACM, New York, 1986,
pp. 50-58. (Erratum in Proceedings of the 20th Annual ACM Symposium on Theory of Computing, 1988,
pp. 551.)
\bibitem{Chertkov}
Michael Chertkov, Lukas Kroc, Massimo Vergassola, {\sl Belief Propagation and Beyond for Particle Tracking}, http://arxiv.org/abs/0806.1199, 2008.
\bibitem{Eg}
G.P. Egorychev, {\sl The solution of van der Waerden's problem for permanents}, Advances
in Math., 42, 299-305, 1981.

\bibitem{Fal}
D. I. Falikman, {\sl Proof of the van der Waerden's conjecture on the permanent of a
doubly stochastic matrix}, Mat. Zametki 29, 6: 931-938, 957, 1981, (in Russian).

\bibitem{FKS}
A. Ferber, M. Krivelevich and B. Sudakov, {\sl Counting and packing Hamilton cycles in dense graphs and oriented graphs}, preprint.

\bibitem{friedland-tverberg}
S. Friedland, {\sl A proof of a generalized van der Waerden conjecture on permanents}, Linear and Multilinear Algebra 11 (1982), no. 2, 107-120.

\bibitem{FRIEDLAND}
S. Friedland, E. Krop, P. H. Lundow, K. Markström, {\sl Validations of the Asymptotic Matching Conjectures}, arxiv preprint arXiv:math/0603001, 2006.

\bibitem{fried}
S. Friedland and L. Gurvits, {\sl Lower Bounds for Partial Matchings in Regular Bipartite Graphs and
Applications to the Monomer-Dimer Entropy}, Combinatorics, Probability and Computing, 2008.

\bibitem{Gurvits-S}
L. Gurvits, {\sl Van der Waerden/Schrijver-Valiant like conjectures and stable (aka hyperbolic) homogeneous polynomials: one theorem for all}, Electronic Journal of Combinatorics 15 (2008).
\bibitem{Gurvits-V}
L. Gurvits, {\sl A polynomial-time algorithm to approximate the mixed volume within a simply exponential factor.}, Discrete Comput. Geom. 41 (2009), no. 4, 533-555.
\bibitem{Jebara}
B. Huang and T. Jebara, {\sl Approximating the Permanent with Belief Propagation.}, New York
Academy of Sciences Machine Learning Symposium 2007. Poster and abstract.
\bibitem{JSV}
M. Jerrum, A. Sinclair, and E. Vigoda, {\sl A polynomial-time approximation algorithm for the permanent of a matrix with nonnegative entries.}, J. ACM 51(4): 671-697 (2004)
\bibitem{Heil}
Heilmann, Ole J.; Lieb, Elliott H. {\sl Theory of monomer-dimer systems.}, Comm. Math. Phys. 25 (1972), 190-232.
\bibitem{LSW}
N. Linial, A. Samorodnitsky, A.Wigderson, {\sl A Deterministic Strongly Polynomial Algorithm for Matrix Scaling and Approximate Permanents.}, Combinatorica 20(4): 545-568 (2000)
\bibitem{L-S} M. Laurent, A. Schrijver, {\sl On Leonid Gurvits's proof for permanents},
American Mathematical Monthly; 117(10):903-911.
\bibitem{Lovasz-Plummer}
L. Lovasz and M. D. Plummer, {\bf Matching Theory}, North Holland, Amsterdam 1986.
\bibitem{MSS}
Adam Marcus, Daniel A. Spielman, Nikhil Srivastava, {\sl Interlacing Families I: Bipartite Ramanujan Graphs of All Degrees},
 arXiv:1304.4132 [math.CO], 2013.
\bibitem{Minc}
H. Minc, {\bf Permanents}, Encyclopeadia of Mathematics and its Applications, vol. 6, Addison-Wesley, Reading, Mass., 1978.

\bibitem{RZ}
M. Rudelson, O. Zeitouni, {\sl Singular values of Gaussian matrices and permanent estimators}, arXiv:1301.6268, 2013.

\bibitem{Sam}
A. Samorodnitsky, {\sl An upper bound for permanents of nonnegative matrices.}, J. Comb. Theory, Ser. A 115(2): 279-292 (2008)


\bibitem{S-V}
A. Schrijver and W.G.Valiant, {\sl On lower bounds for permanents}, Indagationes Mathematicae 42, pp. 425-427, 1980.

\bibitem{Schrijver}
A. Schrijver, {\sl Counting 1-factors in regular bipartite graphs}, Journal of Combinatorial
Theory, Series B 72 (1998) 122-135.

\bibitem{Soules}
G. W. Soules, {\it New permanental upper bounds for nonnegative
  matrices}, Linear and Multilinear Algebra 51, 2003, pp. 319-337.

\bibitem{Valiant}
L. G. Valiant, {\sl The complexity of computing the permanent}, Theoretical Computer
Science, 8(2), 189-201, 1979.

\bibitem{Vontobel}
P.O. Vontobel, {\sl The Bethe permanent of a non-negative matrix}, in Proc. of Commu-
nication, Control, and Computing (Allerton), 2010.


\bibitem{Wigderson-personal}
A. Wigderson, personal communication.

\bibitem{Zygmund}
A. Zygmund, {\bf Trigonometric series}, Volume 1 and 2 combined (3rd ed.), Cambridge University Press, 2002.


\end{thebibliography}
\end{document}